\documentclass[11pt]{article} 

\usepackage[utf8]{inputenc} 

\usepackage{ragged2e}

\usepackage{geometry} 
\usepackage{graphicx} 

\usepackage{booktabs} 
\usepackage{array} 

\usepackage{verbatim} 
\usepackage{subfig} 

\usepackage{fancyhdr} 
\title{The Arithmetic Fourier Transform}
\author{Joel L. Schiff}

\begin{document}

\maketitle

ABSTRACT. The Arithmetic Fourier Transform is a numerical formulation for computing Fourier series and Taylor series coefficients. It competes with the Fast Fourier Transform in terms of speed and efficiency, requiring only addition operations and can be performed by parallel processing. The AFT has some deep connections with the Prime Number Theorem and its rich history is discussed in this expository article.

\bigskip

\begin{center}

\emph{Dedicated to my colleague Wayne J. Walker (1944-2019). Together we uncovered some of the beautiful aspects of the AFT}

\bigskip
\bigskip

\textbf{Table of Contents}

\end{center}

\flushleft 1. Introduction

2. Möbius Function

3. Early Work of Heinrich Bruns

4. Wintner’s Development

5. Analytic Functions

6. Hardy Space 
$H^{1}$ 

7. {$\mu-$}Regular Functions

8. Step-Functions

9. Applications of AFT to Signal Processing

10. References
\bigskip

\justify

\textbf{1. Introduction}
\bigskip

The story of the Arithmetic Fourier Transform is one of discovery, being lost and rediscovered, several times in differing contexts over nearly a century. Thus, its development somewhat parallels the case of the Fast Fourier Transform when it was once again discovered by Cooley and Tukey in 1965 [3]. Indeed, some early work regarding the FFT traces back to Carl Runge in $1903,$ the same year that the AFT saw light of day in the work of Heinrich Bruns [1]. But of course, the great C.F. Gauss was already onto the essence of the FFT by 1805 (cf. [9]). The FFT reworking by Danielson and Lanzcos in 1942 is again in concert with the AFT elaborations in the 1945 monograph of Wintner [24]. Modern extensions of the latter trace to the late 1980s and early 1990s (cf. $[12],[15],[16],[17],[18],[19],[20],[21],[22])$.

One beautiful aspect of the Arithmetic Fourier Transform is that it relates the field of analytic number theory to that of Fourier analysis.
\bigskip

\noindent
\textbf{2. Möbius Function}
\bigskip

Critical to our sampling formula discussion is the Mobius function first investigated by August Ferdinand Möbius in 1832 [11] although it was implicitly known to Euler and touched upon by Gauss (of course). It is defined as,
\bigskip

$\mu(1)=1$;

$\mu(n)=0$ if $p$ is prime such that $p^{2} \mid n$; 

$\mu(n)=(-1)^{j}$ if $n$ is the product of $j$ distinct primes.

\bigskip

\noindent So, for example, $\mu(6)=1$, $\mu(8)=0$, $\mu(30)=-1$.

There is another expression for the Möbius function that does not involve primes (cf. [8]) but rather the \textit{primitive} $n${th} roots of unity (i.e. $z^{n}=1$ and $z^{k} \neq 1$ for any positive integer $k<n$), namely:

\[
\mu(n)=\sum e^{2 \pi i \frac{k}{n}},
\]

\noindent
where the sum is taken over: $1 \leq k \leq n,$ gcd$(k, n)=1 .$ 

A fundamental property of the Möbius function that we will make extensive use of is the following:

\[
\sum_{d \mid n} \mu(d)=\left\{\begin{array}{ll}
1 & \ n=1 \\
0 & \ n>1,
\end{array}\right.
\]
which follows from the preceding expression for $\mu$ since each $n$th root of unity is a primitive root for just one value of $d$ that divides $n$ and the $n$th roots of unity sum to
zero.

The Möbius function also has a strong connection with the Riemann zeta function ${\zeta(s)}$, 
\[
\frac{1}{\zeta(s)}=\sum_{n=1}^{\infty} \frac{\mu(n)}{n^{s}},
\]
for $Re(s)>1.$ Indeed, the Riemann Hypothesis is equivalent to the statement that:
\bigskip
\[
\sum_{n \leq x} \mu(n)=O\left(x^{\frac{1}{2}+\varepsilon}\right),
\]
for all $\varepsilon>0$

What is of interest to us in the present work is the notion of Möbius inversion, a powerful idea that is implicit in the seminal 1832 paper of Möbius and has become a basic tool in the arsenal of number theorists. We state it here in a manner to suit our analytical applications. However, there are other numerous applications of Möbius inversion to problems in physics that can be found in the excellent book by Chen [2].

\bigskip
\noindent
THEOREM 1. \emph{(Möbius inversion). Let $c_{1,}, c_{2}, c_{3}, \cdots,$ be a sequence of complex numbers
such that,
\[
\sum_{n=1}^{\infty}\left|c_{n}\right|<\infty,
\]
and let $b_{1,}, b_{2}, b_{3}, \cdots$ be defined by
\[
b_{n}=\sum_{k=1}^{\infty} c_{k n}.
\]
Under the assumption that,
\[
\sum_{k=1}^{\infty} k\left|c_{k}\right|<\infty,
\]
we have the Möbius inversion formula,}
\begin{equation} 
c_{n}=\sum_{k=1}^{\infty} \mu(k) b_{k n}.
\end{equation}

\noindent
Proof. We consider the array,
\[
\begin{array}{l}
\mu(1) b_{n}= \mu(1) c_{n}+\mu(1) c_{2 n}+\mu(1) c_{3 n}+\mu(1) c_{4 n}+\cdots \\
\mu(2) b_{2 n}=    \quad \hspace{0.6cm}\quad\mu(2) c_{2 n} +\hspace{1.75cm}\mu(2) c_{4 n}+\cdots \\
\mu(3) b_{3 n}=\hspace{2.75cm}\quad \mu(3) c_{3 n}+\cdots \hspace{0.4cm}\quad+\cdots \\
\mu(4) b_{4 n}=\hspace{4.45cm}\quad \mu(4) c_{4 n}+\cdots\\
\quad\vdots
\end{array}
\]
Summing the absolute values of the sum of $j$th column gives:
\[
C_{j}=\left(\sum_{d \mid j}|\mu(j)|\right)\left|c_{j}\right| \leq j\left|c_{jn}\right|.
\]

\noindent
Hence, the sum of the columns satisfies,
\[
\sum_{j=1}^{\infty} C_{j} \leq \sum_{j=1}^{\infty} j\left|c_{jn}\right|\leq\sum_{k=1}^{\infty} k\left|c_{k}\right|<\infty,
\]
and thus by the Weierstrass double series theorem, the double series converges with the sum by columns equal to the sum by rows. As a consequence of the fundamental property of the Möbius function,

\[
\left(\sum_{d \mid j}\mu(d)\right)c_{jn} =0,
\]
for all $j\geq{2}$.

\noindent
Therefore, taking the sum by columns we obtain,
\[
c_{n}=\sum_{k=1}^{\infty} \mu(k) b_{k n},
\]
\noindent
as desired. 

\bigskip

With notation as above, Wintner [24, p.16] presents a more general sufficient condition for Möbius inversion (1), namely if,
\[
\sum_{n=1}^{\infty} 2^{v(n)}\left|c_{n}\right|<\infty
\]
holds, where $v(n)$ is the number of distinct prime divisors of $n$, and more specifically, if
\[
\sum_{n=1}^{\infty} n^{\varepsilon}\left|c_{n}\right|<\infty
\]
for some $\varepsilon>0 .$ He further remarks that the absolute convergence of the terms:

\[
\sum_{n=1}^{\infty}\left|c_{n}\right|<\infty
\]
is in general not sufficient for Möbius inversion, contrary to what was formerly asserted by Bruns.

\bigskip

\noindent
\textbf{3. Early Work of Heinrich Bruns}
\bigskip

The earliest known work in regard to the AFT was by the German mathematician and astronomer Heinrich Bruns [1848-1919] who assumes that a given function $f(x)$ has a Fourier series. Bruns was seeking a more manageable method for calculating the Fourier coefficients in the event that the integral representation for them became impracticable to calculate. Such a situation can occur for instance with a periodic function $f(x)$  that takes values at a \textit{large} number of equidistant points such as daily/hourly astronomical observations of a variable star. Then the harmonic analysis of $f(x)$ will necessarily involve a large partial sum of the Fourier series and the attendant coefficients.

Although Bruns demonstrates Möbius inversion, he does not refer to the Möbius function by its name, only defining it. He then goes on to develop in typical 19th century fashion a formulation (eq. (44) of [1]) for the Fourier cosine coefficients that we can express in a general form as:

\begin{equation}
a_{n}=\sum_{k=1}^{\infty} \frac{\mu(k)}{k n} \sum_{m=1}^{k n} f\left(\frac{m}{k n}\right), \quad n=1,2,3, \ldots
\end{equation}

An expression of the sort on the right-hand side of eq. (2), either in its finite or infinite manifestation, is known as the Arithmetic Fourier Transform. However, Bruns was unable to find a similar expression for the Fourier sine coefficients $b_{n}$ although he developed a cumbersome set of rules for such a purpose.

\bigskip

\noindent
\textbf{4. Wintner's Development}
\bigskip

 The thread was picked up again by Hungarian-American mathematician Aurel Wintner $(1903-1958)$ in his short monograph [24].

We will be interested in taking averages of a function over a suitable domain. Let us start with a function $f: R \rightarrow R$ of period 1 and define,

\begin{equation}
s_{n}(x)=\frac{1}{n} \sum_{m=1}^{n} f\left(x+\frac{m}{n}\right).
\end{equation}
Then $s_{n}(x)$ is also of period 1 and if $f$ is Riemann-integrable,
\[
\lim _{n \rightarrow \infty} s_{n}(x)=\int_{x}^{x+1} f(t) d t=\int_{0}^{1} f(t) d t.
\]
Following Wintner regarding the rate of convergence of the averages $s_{n}(x)$ to the integral we have the following $[24, \mathrm{p.} 4]:$

\bigskip
\noindent
PROPOSITION 1. \emph{Let $f(x)$ be a function of period 1 such that $f^{\prime}(x) \in \ Lip_{1}$ then},

\bigskip

\[
\left|\int_{0}^{1} f(t) d t-\frac{1}{n} \sum_{m=1}^{n} f\left(x+\frac{m}{n}\right)\right| \leq \frac{C}{n^{2}},
\]
\emph{uniformly in $x,$ where $C$ is the Lipschitz constant.}

\bigskip

Recall that $f \in\ Lip_{\alpha}$ on $[a, b]$ with $0<\alpha \leq 1,$ if for some constant $C>0$,
\[
|f(x)-f(y)| \leq C|x-y|^{\alpha}
\]
for all $x, y \in[a, b].$

Wintner subsequently established the following result  [24,\ p.6 ] connecting the Möbius function with Fourier series based upon his estimate of Proposition $1 .$

\bigskip
\noindent
THEOREM $2 .$ \emph{Let $f: R \rightarrow R$ be a function of period 1 and $f^{\prime}(x) \in \ Lip_{1}$ with the normalization,
\[
a_{0}=\int_{x}^{x+1} f(t) d t=\int_{0}^{1} f(t) d t=0,
\]
and let $s_{n}(x)$ be defined as in $(3) .$ Then each of the series,}

\begin{equation}
S_{n}(x)=\sum_{k=1}^{\infty} \mu(k) s_{k n}=\sum_{k=1}^{\infty} \frac{\mu(k)}{k n} \sum_{m=1}^{k n} f\left(x+\frac{m}{k n}\right), \quad n=1,2,3, \ldots
\end{equation}
\noindent
\emph{is absolutely uniformly comvergent, and
\[
f(x)=\sum_{n=1}^{\infty} S_{n}(x)
\]
is absolutely uniformly convergent, with,}

\begin{equation}
S_{n}(x)=a_{n} \cos 2 \pi n x+b_{n} \sin 2 \pi n x.
\end{equation}

Setting $x=0$, in both (4) and $(5),$ we arrive at the result of Bruns, namely,

\begin{equation}
a_{n}=\sum_{k=1}^{\infty} \frac{\mu(k)}{k n} \sum_{m=1}^{k n} f\left(\frac{m}{k n}\right), \quad n=1,2,3, \ldots
\end{equation}

\noindent Setting $x=\frac{1}{4 n}$ in $(5),$ we do obtain another formula for each $b_{n}$ on another set of points, but again no fixed $x$ (or finite fixed set of $x-$values) will supply the $b_{n}$ coefficients for all $n.$

There is a remarkable result of Wintner that is similar but not quite the AFT that is intimately connected with the Prime Number Theorem. Let us first consider the following case.

\bigskip
\noindent
THEOREM $3 .$ \emph{Let $f(x)$ be a Riemann-integrable function of period $1 .$ Then,}

\begin{equation}
\int_{0}^{1} f(t) d t=\sum_{n=1}^{\infty} \frac{1}{n} \sum_{d \mid n} \mu(d) f\left(\frac{n x}{d}\right),
\end{equation}
\emph{for every irrational number $x$.}
\bigskip

To appreciate the deep connection with the PNT, consider the period 1 function,
\[
f(t)=\left\{\begin{array}{cc}
1 & t=x \\
0 & 0<t<x ; \quad x<t<1.
\end{array}\right.
\]
The value at the endpoints $f(0)=f(1)$ is of no consequence as they do not figure into the sum on the right-hand side of eq.(7). Note that the interior sum runs over all the divisors $d$ of $n$ so that $f\left(\frac{n x}{d}\right)=0$ except for $d=n$ since $x$ is irrational. As $f(x)=1,$ we have by (7),
\[
0=\sum_{n=1}^{\infty} \frac{\mu(n)}{n},
\]
and the convergence of $\sum_{n=1}^{\infty} \mu(n) / n$ is a well-known equivalent of the PNT (cf. the monograph of Landau [10] where various equivalents of the PNT are discussed).

The formulation in eq.(7) leads directly to the more general form (Wintner [24], p.24) for Fourier series.

\bigskip

\noindent
THEOREM 4. \emph{Let $f(x)$ be a Riemann-integrable function of period 1 with Fourier coefficients $c_{k}, k=0,\pm 1,\pm 2, \ldots$ Then},

\begin{equation}
c_{k}=\sum_{n=1}^{\infty} \frac{1}{n} \sum_{d \mid n} \mu(d) f\left(\frac{n x}{d}\right) e^{-\frac{2 \pi i k n x}{d}},
\end{equation}

\noindent\emph{where $x$ is an arbitrary irrational mumber.}

\bigskip

The PNT will appear again in the sequel regarding a formula due to Davenport.
Before we apply Möbius inversion to other classes of functions, we require:

\bigskip
\noindent
LEMMA 1.\emph{For any positive integer $j,$}
\[
\sum_{m=1}^{n} e^{i j\left(\frac{2 \pi m}{n}\right)}=\left\{\begin{array}{ll}
0 &  j \neq k n \\
n & \ j=k n,
\end{array}\right.
\]
\emph{for some positive integer $k .$}

\bigskip

The proof is elementary and will be omitted.

\bigskip

\noindent
\textbf{5. Analytic Functions}
\bigskip

We now consider the case of an analytic function $f(z)$ on the closed unit disk $\bar{U}$, hence in a slightly larger open disk containing $U$. It will be demonstrated that each of the Taylor coefficients of $f(z)$ can be determined by an infinite series of the values of the function $f(z)$ taken at symmetrically arrayed sets of points on the boundary $\partial U:|z|=1 .$ The AFT formula (9) was first given for the Taylor coefficients of an analytic function in Schiff -Walker [15], completely unaware of the existence of the work of Bruns or Wintner. The proof there was reliant on an obscure result of J.L. Walsh [23] but the proof given below follows the much more refined version given in [18]. Again, in the sequel, without loss of generality, we will take the normalization,

\[
c_{0}=\frac{1}{2 \pi} \int_{0}^{2 \pi} f\left(e^{i \theta}\right)=0.
\]

\noindent
THEOREM $5.$ \emph{If $f(z)$ is analytic on the closed unit disk $\bar{U}$, with Taylor series,}

\[
f(z)=\sum_{j=1}^{\infty} c_{j} z^{j},
\]
\emph{then each Taylor coefficient is given by,}

\begin{equation}
c_{n}=\sum_{k=1}^{\infty} \frac{\mu(k)}{k n} \sum_{m=1}^{k n} f\left(e^{i\left(\frac{2 \pi m}{k n}\right)}\right), \quad n=1,2,3, \ldots
\end{equation}

\noindent Proof. On $|z|=1$ we can write,
\[
f(\theta)=f\left(e^{i \theta}\right)=\sum_{j=1}^{\infty} c_{j} e^{i j \theta}.
\]
Define the averages over a symmetrically arrayed set of points,
\[
f_{n}(\theta)=\frac{1}{n} \sum_{m=1}^{n} f\left(\theta+\frac{2 \pi m}{n}\right), \quad n=1,2,3, \ldots
\]
Then,
\[
f_{n}(\theta)=\frac{1}{n} \sum_{m=1}^{n} \sum_{j=1}^{\infty} c_{j} e^{i j\left(\theta+\frac{2 \pi m}{n}\right)}
\]

\[
\begin{array}{l}
\bigskip
=\displaystyle\sum_{j=1}^{\infty} c_{j} e^{i j \theta} \frac{1}{n} \sum_{m=1}^{n} e^{i j\left(\frac{2 \pi m}{n}\right)} \\
=\displaystyle\sum_{k=1}^{\infty} c_{k n} e^{i kn \theta},
\end{array}
\]
{by Lemma 1.}

Upon setting $\theta=0$, we now have,
\[
f_{n}=f_{n}(0)=\sum_{k=1}^{\infty} c_{k n}.
\]
By the analyticity of $f(z)$,
\[
\sum_{k=1}^{\infty} k\left|c_{k}\right|<\infty,
\]
and hence by Theorem $1,$ we have,

\[
\begin{array}{c}
\bigskip
\displaystyle c_{n}=\sum_{k=1}^{\infty} \mu(k) f_{k n} \\
\displaystyle=\sum_{k=1}^{\infty} \frac{\mu(k)}{k n} \sum_{m=1}^{k n} f\left(e^{i\left(\frac{2 \pi m}{k n}\right)}\right),
\end{array}
\]
as desired.
\bigskip

\noindent
COROLLARY$1 .$ \emph{If $f(z)$ is analytic on the open unit disk $U:[z]<1,$ with Taylor series,}
\[
f(z)=\sum_{n=1}^{\infty} c_{n} z^{n},
\]
\emph{then the Taylor coefficients are given by,}

\begin{equation}
c_{n}=\frac{1}{r^{n}} \sum_{k=1}^{\infty} \frac{\mu(k)}{k n} \sum_{m=1}^{k n} f\left(r e^{i\left(\frac{2 \pi m}{k n}\right)}\right), \quad 0<r<1, \quad n=1,2,3, \ldots
\end{equation}
In the case of the Z-transform given by $X(z)=\sum_{j=1}^{\infty} c_{j} z^{-j},$ we have by considering the function $w=\frac{1}{z}$:

\bigskip
\noindent
COROLLARY $2 .$ \emph{If}
\[
X(z)=\sum_{j=1}^{\infty} c_{j} z^{-j}
\]
\emph{converges for $|z|>r$ and $r<1,$ then}

\begin{equation}
c_{n}=\sum_{k=1}^{\infty} \frac{\mu(k)}{k n} \sum_{m=1}^{k n} X\left(e^{-i\left(\frac{2 \pi m}{k n}\right)}\right), \quad n=1,2,3, \ldots
\end{equation}
Note that the AFT only involves addition/subtraction operations so that it can be efficiently computed by parallel processing.

\bigskip

\noindent
\textbf{6. Hardy Space} \boldmath$H^{1}$ 
\unboldmath
\bigskip

We now consider the Hardy space $H^{1}(U)$, that is, functions $f(z)$ that are analytic in the unit disk $U:|z|<1$ and satisfy the condition,

\[
\sup _{0 \leq r<1} \int_{0}^{2 \pi}\left|f\left(r e^{i \theta})\mid d \theta<\infty\right.\right..
\]

Note that by changing the Wintner setting to a real-valued function $f(z)$ of period $2\pi$ defined on the circle $|z| = 1$ such that $f^{\prime} \in L i p_{1},$ we can re-write the Fourier cosine coefficients in (6) as,

\begin{equation}
a_{n}=\sum_{k=1}^{\infty} \frac{\mu(k)}{k n} \sum_{m=1}^{k n} f\left(e^{i\left(\frac{2 \pi m}{k n}\right)}\right)\quad n=1,2,3, \ldots
\end{equation}
where again we are taking the integral mean value to be zero.

\bigskip

\noindent
THEOREM $6 .$ \emph{ If $f(z) \in H^{1}(U)$ and $f^{\prime}\in\ Lip_{1}$ on $|z|=1,$ then the Taylor/Fourier coefficients are given by}
                                                                
\begin{equation}
c_{n}=\sum_{k=1}^{\infty} \frac{\mu(k)}{k n} \sum_{m=1}^{k n} f\left(e^{i\left(\frac{2 \pi m}{k n}\right)}\right), \quad n=1,2,3, \ldots
\end{equation}

\noindent Proof. From the hypothesis, we can infer that $f(\theta)=f\left(e^{i \theta}\right)$ has a Fourier series expansion at each point of $|z|=1$,
\[
\displaystyle{f(\theta)=\sum_{n=1}^{\infty} c_{n} e^{i n \theta}},
\]
and $c_{-n}=0$ for $n=1,2,3, \ldots .$ with the $c_{n}$'s also the Taylor series coefficients (cf. $[6,$
p.38]). Therefore,
\[
\begin{array}{c}

\displaystyle{f(\theta)=\sum_{n=1}^{\infty}\left(a_{n}+i b_{n}\right) e^{i n \theta}} \\
\bigskip

=\displaystyle{\sum_{n=1}^{\infty}\left(a_{n} \cos n \theta-b_{n} \sin n \theta\right)+i \sum_{n=1}^{\infty}\left(b_{n} \cos n \theta+a_{n} \sin n \theta\right)} \\

u\left(e^{i \theta}\right)+i v\left(e^{i \theta}\right).
\end{array}
\]

Here, the functions $u, v$ satisfy the Lipschitz condition on their derivatives which by (6) yields,
\[
a_{n}=\sum_{k=1}^{\infty} \frac{\mu(k)}{k n} \sum_{m=1}^{k n} u\left(e^{i\left(\frac{2 \pi m}{k n}\right)}\right), \quad b_{n}=\sum_{k=1}^{\infty} \frac{\mu(k)}{k n} \sum_{m=1}^{k n} v\left(e^{i\left(\frac{2 \pi m}{k n}\right)}\right)
\]
for $n=1,2,3, \ldots .$ Then,
\[
c_{n}=a_{n}+i b_{n}=\sum_{k=1}^{\infty} \frac{\mu(k)}{k n} \sum_{m=1}^{k n} f\left(e^{i\left(\frac{2 \pi m}{k n}\right)}\right), \quad n=1,2,3, \ldots
\]
proving the theorem.

\bigskip
An error estimate is given by:

\bigskip
\noindent
COROLLARY 3. \emph{If $f(z) \in H^{1}(U)$ and $f^{\prime}\in\ Lip_{1}$ on $|z|=1,$ with Lipschitz constant
C. Then the truncation error for each Taylor/Fourier coefficient is,}
\bigskip
\[
\left|c_{n}-\sum_{k=1}^{N} \frac{\mu(k)}{k n} \sum_{m=1}^{k n} f\left(e^{i\left(\frac{2 \pi m}{k n}\right)}\right)\right| \leq \frac{C}{n^{2} N} .
\]

\bigskip

In fact by Proposition 1,
\[
\begin{array}{l}
\bigskip
\displaystyle\left|c_{n}-\sum_{k=1}^{N} \frac{\mu(k)}{k n} \sum_{m=1}^{k n} f\left(e^{i\left(\frac{2 \pi m}{k n}\right)}\right)\right| \\
\displaystyle\quad \leq \frac{C}{n^{2}} \sum_{k=N+1}^{\infty} \frac{1}{k^{2}} \leq \frac{C}{n^{2} N}.
\end{array}
\]
In practice, computer simulations of the finite AFT have shown a much more rapid convergence than this (cf. $[12],[18]) .$

\bigskip

\indent Theorems 5 and 6 have underlying different assumptions both bestowing the AFT. While Theorem 5 requires the function $f(z)$ to be analytic in a slightly larger disk containing $U:|z|<1,$ Theorem 6 requires an integrability condition in $U$ and a
Lipschitz condition on $\partial U$.

Indeed, it is possible to derive Theorem 6 using the approach of Theorem $5,$ specifically as formulated in Corollary 1. We utilize a result of Hardy and Littlewood, namely that if $f^{\prime} \in Lip_{1}$ on $|z|=1,$ then there exists a constant $C$ for $0<r \leq 1,$ with (cf. $[7, p .413])$
\[
\left|f^{\prime}\left(r e^{i \theta}\right)-f^{\prime}\left(r e^{i \theta^{\prime}}\right)\right| \leq C\left|r e^{i \theta}-r e^{i \theta^{\prime}}\right|.
\]

Applying the error estimate of Proposition 1 to each circle of radius $r$, we obtain,

\begin{equation}
\left|\frac{1}{k n} \sum_{m=1}^{k n} f\left(r e^{i\left(\frac{2 \pi m}{k n}\right)}\right)\right| \leq \frac{C}{k^{2} n^{2}}, \quad 0<r \leq 1.
\end{equation}

Let $\varepsilon>0 .$ From eq.(10), for a fixed radius $r_{0}$ and some $N$ sufficiently large such that for $r \geq r_{0}$, applying (14) as in Corollary 3,
\[
\left|c_{n}-\sum_{k=1}^{\infty} \frac{\mu(k)}{k n} \sum_{m=1}^{k n} f\left(e^{i\left(\frac{2 \pi m}{k n}\right)}\right)\right|
\]

\[
\begin{array}{l}

\bigskip
\displaystyle=\left|\frac{1}{r^{n}} \sum_{k=1}^{\infty} \frac{\mu(k)}{k n} \sum_{m=1}^{k n} f\left(r e^{i\left(\frac{2 \pi m}{k n}\right)}\right)-\sum_{k=1}^{\infty} \frac{\mu(k)}{k n} \sum_{m=1}^{k n} f\left(e^{i\left(\frac{2 \pi m}{k n}\right)}\right)\right| \\

\displaystyle<\frac{\varepsilon}{2}+\left|\frac{1}{r^{n}} \sum_{k=1}^{N} \frac{\mu(k)}{k n} \sum_{m=1}^{k n} f\left(r e^{i\left(\frac{2 \pi m}{k n}\right)}\right)-\sum_{k=1}^{N} \frac{\mu(k)}{k n} \sum_{m=1}^{k n} f\left(e^{i\left(\frac{2 \pi m}{k n}\right)}\right)\right|
\end{array}
\]
Now the second term involves only a finite number of terms and so it can be made
less than $\varepsilon / 2$ for $r$ sufficiently close to $1$. It follows that,
\[
\left|c_{n}-\sum_{k=1}^{\infty} \frac{\mu(k)}{k n} \sum_{m=1}^{k n} f\left(e^{i\left(\frac{2 \pi m}{k n}\right)}\right)\right|<\varepsilon,
\]
as desired.

\bigskip
\noindent
 Remark. For $f^{\prime}\in \ Lip_{1}$ on $|z|=1,$ where the Fourier series coefficients satisfy $c_{n}=0$ for $n<-k, k>0,$ we may consider the function $g(z)=z^{k} f(z)$ and obtain similar formulae for the $c_{n}$'s, $n \geq-k$.

\bigskip
\noindent
\boldmath
{7. $\mu-$}\textbf{Regular Functions}
\unboldmath

\bigskip
A potential theory was developed by R.J. Duffin [5] pertaining to the Yukawan potential for the strong nuclear force, $\frac{e^{-\mu r}}{r}$. Such potentials satisfy the elliptic p.d.e.

\begin{equation}
\Delta u=\mu^{2} u \quad(\mu>0)
\end{equation}

\noindent and any $C^{2}-$ solution of (15) in the complex plane is called panharmonic. We form the pseudoanalytic function $f=u+i v$ with $u, v \in C^{2}$ satisfying the following \lq{Cauchy-Riemann equations}\rq,

\[\displaystyle\frac{\partial u}{\partial x}=\frac{\partial v}{\partial y}+\mu u\]

\[\displaystyle\frac{\partial u}{\partial y}=-\frac{\partial v}{\partial x}-\mu v.\]

Then $u, v$ are \emph{panharmonic} and $f(z)$ is called \emph{$\mu$-regular.} Note that as $\mu \rightarrow 0,$ a $\mu-$regular function becomes analytic in the classical sense. According to Duffin $[5], \mu-$regular functions have the following Fourier series analogue.

\bigskip
\noindent
THEOREM $7.$ \emph{Let $f(z)$ be $\mu-$regular in the closed unit disk $\bar{U}$ (that is, in a slightly larger open disk containing $\bar{U}$). Then for $|z|=r \leq 1$,}
\[
f(z)=\sum_{n=-\infty}^{\infty} c_{n} I_{|n|}(\mu r) e^{i n \theta}
\]
\emph{where,
\[
c_{n}=\frac{1}{2 \pi I_{n}(\mu)} \int_{0}^{2 \pi} f\left(e^{i \theta}\right) e^{-i n \theta} d{\theta}, \quad n=0,1,2,3, \ldots
\]
and $I_{n}$ is the modified Bessel function of the first kind,}
\[
I_{n}(x)=\frac{1}{n !}\left(\frac{x}{2}\right)^{n}\left[1+\frac{(x / 2)^{2}}{1 \cdot(n+1)}+\frac{(x / 2)^{4}}{1 \cdot 2 \cdot(n+1)(n+2)}+\cdots\right].
\]
Furthermore, the Fourier coefficients having negative index satisfy the condition:
\[
c_{-n}=\overline{c_{n-1}},
\]
\noindent
a characterizing feature of $\mu-$regular functions.
Such functions also satisfy a mean-value property, i.e. for $n=0$ and $0 \leq \rho<1$,
\[
f(0)=\frac{1}{2 \pi I_{0}(\mu \rho)} \int_{0}^{2 \pi} f\left(\rho e^{i \theta}\right) d \theta.
\]
If $f(z)$ is as in Theorem $7,$ and without loss of generality $f(0)=0$, then the AFT for $\mu-$regular functions is given by the recursive formula [17] :

\begin{equation}
c_{n}=\frac{1}{I_{n}(\mu)} \sum_{k=1}^{\infty} \frac{\mu(k)}{k n} \sum_{m=1}^{k n}f\left(e^{\frac{2 \pi i m}{k n}}\right)-\overline{c_{n-1}}
\end{equation}
$n=1,2,3, \ldots$

\bigskip
\noindent
\textbf{8. Step-Functions}

\bigskip

Although the AFT would seem to require a high degree of regularity, remarkably the canonical representation can be demonstrated for even step-functions. These are important in the field of signal processing and time-series analysis. Once again, a deep result from number theory due to Harold Davenport will come into play.

The number theoretical formula in question is based on the \emph{first Bernoullian function},

\begin{equation}
\{t\}=\left\{\begin{array}{ll}
\bigskip
t-[t]-\frac{1}{2}, & t \neq[t] \\
0, & t=[t],
\end{array}\right.
\end{equation}

\noindent where $[t]$ is the greatest integer part of $t .$ This function was employed by Davenport to establish the \emph{Davenport formula},

\begin{equation}
\sum_{k=1}^{\infty} \frac{\mu(k)}{k}\{k \theta\}=-\frac{1}{\pi} \sin 2 \pi \theta,
\end{equation}

\noindent with the convergence being uniform in $\theta$. Davenport established eq. (18) using the Prime Number Theorem and some work of Vinogradov. Conversely, it can be shown that (18) implies the PNT [19].

Indeed, for any $\varepsilon>0$, there is a positive integer $N=N(\varepsilon)$ sufficiently large such that by (18)
\[
\left|-\frac{1}{\pi} \sin 2 \pi \theta-\sum_{k=1}^{N} \frac{\mu(k)}{k}\{k \theta\}\right|<\varepsilon,
\]
for all $\theta \in[0,2 \pi]$. Furthermore, for all sufficiently small $\theta>0$, we can make $\left|-\frac{1}{\pi} \sin 2 \pi \theta\right|<\varepsilon,$ implying that
\[
\left|\sum_{k=1}^{N} \frac{\mu(k)}{k}\{k \theta\}\right|<2 \varepsilon.
\]
Choosing $\theta>0$ possibly smaller, we can obtain,
\[
\{k \theta\}=k \theta-\frac{1}{2}, \quad k=1,2, \ldots, N,
\]
and $\theta<\frac{\varepsilon}{N} .$ It follows that for $\theta$ sufficiently small,
\[
\begin{array}{c}
\displaystyle 2 \varepsilon>\left|\sum_{k=1}^{N} \frac{\mu(k)}{k}\{k \theta\}\right|=\left|\sum_{k=1}^{N} \frac{\mu(k)}{k} k \theta-\frac{1}{2} \sum_{k=1}^{N} \frac{\mu(k)}{k}\right| \\ \\

\displaystyle\geq \frac{1}{2}\left|\sum_{k=1}^{N} \frac{\mu(k)}{k}\right|-\left|\sum_{k=1}^{N} \mu(k) \theta\right|.
\end{array}
\]

\noindent
As a consequence,

\[
\begin{array}{c}
\displaystyle\frac{1}{2}\left|\sum_{k=1}^{N} \frac{\mu(k)}{k}\right|<2\varepsilon+\sum_{k=1}^{N} |\mu(k)| \theta
\end{array}
\]

\[
<2 \varepsilon+\theta N<3 \varepsilon.
\]

\noindent Therefore,
\[
\sum_{n=1}^{\infty} \frac{\mu(k)}{k}=0,
\]
and as mentioned in Sec. 4, this implies the PNT.

A further connection between the Prime Number Theorem via the Davenport formula and sampling theory is provided by the following.

\bigskip
\noindent
THEOREM $8.$ \emph{Let $f(\theta)$ be an even step-function defined on $[-\pi, \pi]$ and extended to be periodic of period $2 \pi .$}

\bigskip

\emph{(i) Suppose that $f(\theta)$ is normalized so that,}
\[
\int_{-\pi}^{\pi} f(\theta) d \theta=0.
\]

\emph{(ii) At a discontinuity $\theta$,
\[
f(\theta)=\frac{1}{2}(f(\theta+)-f(\theta-)).
\]
Then the Fourier cosine coefficients of $f(\theta)$ are given by,}
\[
a_{n}=\sum_{k=1}^{\infty} \frac{\mu(k)}{k n} \sum_{m=0}^{k n-1} f\left(\frac{2 \pi m}{k n}\right), \quad n=1,2,3, \ldots
\]
Proof. For simplicity let,

\begin{equation}
F_{N}=\sum_{m=0}^{N-1} f\left(\frac{2 \pi m}{N}\right).
\end{equation}

\noindent First, we consider the specific even step-function,

\[
f_{b}(\theta)=\left\{\begin{array}{ll}

\bigskip
1-\frac{b}{\pi}, & |\theta|<b \\

\bigskip

\frac{1}{2}-\frac{b}{\pi}, & \theta=\pm b \\
-\frac{b}{\pi}, & b<\theta \leq \pi.
\end{array}\right.
\]

\noindent This function is simply the normalization of the function,
\[
g_{b}(\theta)=\left\{\begin{array}{ll}

\bigskip
1, & |\theta|<b \\

\bigskip
\frac{1}{2}, & \theta=\pm b \\
0, & b<\theta \leq \pi,
\end{array}\right.
\]
by subtracting the constant,
\[
a_{0}=\frac{1}{2 \pi} \int_{-\pi}^{\pi} g_{b}(\theta) d \theta=\frac{b}{\pi}.
\]
For $n \geq 1$, will have the same Fourier cosine coefficients for both $f_{b}$, and $g_{b}$, namely,
\[
\begin{array}{c}
\bigskip
\displaystyle{a_{n}=\frac{1}{\pi} \int_{-\pi}^{\pi} \cos n \theta g_{b}(\theta) d \theta} \\
\displaystyle{=\frac{2}{n \pi} \sin n b}.
\end{array}
\]
In the special case of the function $f_{b}(\theta)$ above, we must therefore show that,

\begin{equation}
\frac{2}{n \pi} \sin n b=\sum_{k=1}^{\infty} \frac{\mu(k)}{k n} F_{k n}
\end{equation}

\noindent for $F_{N}$ given by (19).
We claim that $F_{N}=-2\left\{\frac{b N}{2 \pi}\right\} .$ To see this, we consider two cases.

\bigskip

\noindent Case $1 .$ Suppose that $b=\frac{2 \pi m}{N}$ for some positive integer $m$ so that $b$ is itself a sample point and therefore,
\[
\left\{\frac{b N}{2 \pi}\right\}=\{m\}=0.
\]
Furthermore, by conditions (i) and (ii), a simple calculation shows that,
\[
F_{N}=\sum_{m=0}^{N-1} f_{b}\left(\frac{2 \pi m}{N}\right)=\frac{N}{2 \pi} \int_{0}^{2 \pi} f_{b}(\theta) d \theta=0.
\]
Case $2 .$ Next, we assume that $b$ is not a sample point. In this case, the number of positive multiples of $2 \pi / N$ that are less than $b$, i.e. sample points, is $\left[\frac{b}{2 \pi} N\right] .$ Hence, on $(-b, b)$, including the origin, there are $2\left[\frac{b}{2 \pi} N\right]+1$ sample points where $f_{b}$ takes the value $\left(1-\frac{b}{\pi}\right),$ as well as $N-2\left[\frac{b}{2 \pi} N\right]-1$ points where $f_{b}$ takes the value $\left(-\frac{b}{\pi}\right)$. As a consequence, invoking (17) at the last step,
\[
\begin{array}{c}
\bigskip
F_{N}=\left(2\left[\frac{b N}{2 \pi}\right]+1\right)\left(1-\frac{b}{\pi}\right)+\left(N-2\left[\frac{b N}{2 \pi}\right]-1\right)\left(-\frac{b}{\pi}\right) \\

\bigskip
=2\left[\frac{b N}{2 \pi}\right]-\frac{b N}{\pi}+1 \\

=-2\left\{\frac{b N}{2 \pi}\right\},
\end{array}
\]
as claimed.

Hence by the Davenport formula (18),
\[
\begin{array}{c}
\bigskip

\displaystyle\sum_{k=1}^{\infty} \frac{\mu(k)}{k n} F_{k n}=-2 \sum_{k=1}^{\infty} \frac{\mu(k)}{k n}\left\{\frac{knb}{2 \pi}\right\} \\
\displaystyle=\frac{2}{n \pi} \sin nb,
\end{array}
\]
as required for (20). This proves the special case for the given even step-function $f_{b}(\theta)$.
More generally, every even step-function will be a finite linear combination of the constant function and functions of the form $g_{b}(\theta)$. The normalization of the constant function is zero, and $f_{b}(\theta)$ is the normalization of $g_{b}(\theta)$. As a consequence, it suffices to consider a linear combination of step-functions $f_{b}(\theta),$ which now establishes the theorem in the general case.

\bigskip
\noindent
\textbf{9. Application of AFT to Signal Processing}

\bigskip

From a signal processing point of view a similar expression for the Fourier cosine coefficients as in eq.(2) appears in for example Tufts and Sadasiv [21]. Indeed, let $F(t)$ be a real-valued function of period $1,$ given by the Fourier series

\[
F(t)=\sum_{k=1}^{\infty} a_{n}(t)
\]
where,
\[
a_{n}(t)=A_{n} \cos \left(2 \pi n t+\theta_{n}\right).
\]
Here $A_{n}$ and $\theta_{n}$ are the amplitude and phase respectively of the $n$th harmonic and as previously we assume that $F(t)$ has no constant term. Moreover, it is assumed that all harmonic terms above the $N^{th}$ term are zero, that is,

\[
a_{n}=0,  \quad n=N+1,N+2, \ldots
\]

Furthermore, define a bank of $N$ delay-line filters each having input $F(t).$ Denoting the output of the $k$th filter by $s_{k}(t)$, we define as in (3),
\[
s_{k}(t)=\frac{1}{k} \sum_{m=0}^{k-1} F\left(t-\frac{m}{k}\right),
\]
\noindent
for $k=1,2, \ldots N,$ and

\begin{equation}
s_{k}=0,  \quad n=N+1,N+2, \ldots
\end{equation}

\noindent Using Möbius inversion and Lemma $1,$ the authors show that,
\[
\begin{array}{c}
\bigskip

\displaystyle a_{n}(t)=\sum_{k=1}^{\infty} \mu(k) s_{k n}(t) \\

\displaystyle=\sum_{k=1}^{\infty} \frac{\mu(k)}{k n} \sum_{m=0}^{k n-1} F\left(t-\frac{m}{k n}\right), \quad n=1,2, \ldots N.
\end{array}\]

\noindent The sum here is of course finite containing at most $N$ terms due to the assumption (21).

The authors further present an instructive VLSI implementation noting the AFT's computational advantages over the Fast Fourier Transform. Further extensions of Wintner's main result to signal processing are found in Reed \textit{et  al.} $[12] .$ In the words of these authors the AFT ``competes with the classical FFT approach in terms of accuracy, complexity, and speed."

\bigskip
\noindent
\textbf{10. References}

\bigskip

1. H. Bruns, \emph{Grundlinien des wissenschaftlichen Rechnens}, Leipzig, 1903.

2. N. Chen,  \emph{Möbius Inversion in Physics, World Scientific}, 2010.

3. J.W. Cooley and J.W. Tukey,  \emph{An algorithm for the machine calculation of complex Fourier series}, Math. Comp., $19(1965), 297-301$.

4. H. Davenport,  \emph{On some infinite series involving arithmetical functions}, I and II, Quart. J. Math. $8(1937), 8-13$ and $313-320$.

5. R.J. Duffin,  \emph{Yukawan potential theory}, J. Math. Anal. Appl. 35 (1971), 104-130.

6. P.L. Duren,  \emph{Theory of $H^{p}$ spaces}, Academic Press, 1970.

7. G.M. Goluzin,  \emph{Geometric Theory of Functions of a Complex Variable,} Amer. Math. Soc. Providence, $\mathrm{RI}, 1969.$

8. G.H. Hardy and E.M. Wright,  \emph{An Introduction to the Theory of Numbers}, Oxford,
1980.

9. M.T. Heideman, D.H. Johnson, and C.S. Burrus,  \emph{Gauss and the history of the Fast Fourier Transform}, IEEE ASSP Magazine, $\mathbf{1}(4)(1984), 14-21.$

10. E. Landau,  \emph{Handbuch der Lehre von der Verteilung der Primzahlen} Teubner, Leipzig, 1909; reprinted (with an appendix by P. T. Bateman), Chelsea, New York, 1953.

11. A.F. Möbius,  \emph{Über eine besondere Art von Umkehrung der Riehen}, J. reine angewandte Math., $9(1832), 105-123$

12. I.S. Reed, D.W. Tufts, X. Yu, T.K. Truong, M.-T Shih, and X. Yin,  \emph{Fourier analysis and signal processing by use of the Möbius imversion formula}, IEEE Trans. ASSP $38(1990), 458-470$.

13. C. Runge,  \emph{Über die zerlegung empirisch gegebener periodischer funktionen in simuswellen}, Zeit. Math. Phys., $49,443-456,1903$.

14. C. Runge,  \emph{Über die zerlegung einer empirischen funktion in simuswellen}, Zeit. Math. Phys., $53,117-123,1905$.

15. J.L. Schiff and W.J. Walker,  \emph{A sampling theorem for analytic functions}, Proc. Amer. Math. Soc., $99(1987), 737-740$.

16. J.L. Schiff and W.J. Walker,  \emph{A sampling theorem and Wintner 's results on Fourier coefficients}, J. Math. Anal. Appl., 133(1988), 

466-471.

17. J.L. Schiff and W.J. Walker,  \emph{A sampling theorem for a class of pseudoanalytic functions}, Proc. Amer. Math. Soc. $111(1991), 695-699 .$

18. J.L. Schiff, T.J. Surendonk, and W.J. Walker,  \emph{An algorithm for computing the inverse Z Transform}, IEEE, $40(1992), 2194-2198$.

19. J.L. Schiff and W.J Walker,  \emph{A sampling formula in signal processing and the prime number theorem,} NZ J. Math., $23(1994) 147-155$.

20. D.W. Tufts,  \emph{Comments on A note of the computational complexity of the Arithmetic Fourier Transform}, IEEE Trans. ASSP, 37 1147-1148, 1989.

21. D.W. Tufts and G. Sadasiv, \emph{The Arithmetic Fourier Transform}, IEEE ASSP Magazine, $5(1988), 13-17$.

22. D.W. Tufts, Z. Fan, and Z. Cao, \emph{ Image processing and the Arithmetic Fourier Transform,} SPIE/IST Conference, Jan. 1989 .

23. J.L. Walsh,  \emph{A mean value theorem for polynomials and harmonic polynomials,} Bull. Amer. Math. Soc. $42(1936), 923-930$.

24. A. Wintner,  \emph{An Arithmetical Approach to Ordinary Fourier Series}, Waverly Press, 1945.

\bigskip
\bigskip

\emph{AMS Mathematics Subject Classification:} Primary 42A16, 11A25, Secondary 30J99.

\bigskip
Department of Mathematics

University of Auckland

Auckland, New Zealand

\end{document}